\documentclass[10pt, a4paper]{article}

\usepackage{amsmath,  amsfonts}

\usepackage[T1]{fontenc}

\newcommand{\be} {\begin{equation}}
\newcommand{\ee} {\end{equation}}

\begin{document}

\title{An inequality of type sup+inf on $ {\mathbb S}_4 $ for the Paneitz operator.}
\author{Samy Skander Bahoura \\
\small { Equipe d'Analyse Complexe et G\'eom\'etrie } \\
\small { Universit\'e Pierre et Marie Curie, 4 Place Jussieu, 75005, Paris, France} \\
\small {e-mails: samybahoura@yahoo.fr, samybahoura@gmail.com}
}
\date{}

\maketitle

\begin{abstract} We give an inequality of type $ \sup+\inf $ on ${\mathbb S}_4 $ for the Paneitz equation.
\end{abstract}

\section{Introduction and Main Result.} 

We set $ \Delta = - \nabla^i(\nabla_i) $ the Laplace-Beltrami operator.

\bigskip

On $ {\mathbb S}_4 $ we consider the following Paneitz equation:

$$ \Delta^2 u +2 \Delta u + 6 = V e^{4u} \qquad (E)  $$

with, $ 0 \leq V(x) \leq b, \,\, x\in {\mathbb S}_4 $ and $ b > 0 $.

\bigskip

Here, we try to derive a minoration of the sum $ \sup + \inf $ with minimal condition on $ V $.

\bigskip

Not that, Paneitz (see [7]) had introduced this equation which is invariant under conformal change of metrics.

\bigskip

This equation is similar to the Gaussian curvature equation in dimension 2, which is:

$$ \Delta u+ K =\tilde K e^{2u} \qquad (E') $$

where $ K $ is the Gaussian curvature and $ \tilde K=V/2 $ is the prescribed Gaussian curvature. ($ V $ is the prescribed scalar curvature).

\bigskip

We refer to [2,7], for some examples where this equation where considered and studied.

\bigskip

Some authors where interested by this equation and established  some Sobolev inequalities, see [2-5,7].

\bigskip

For the corresponding equation in the heigher dimensionnal case to $ (E) $ (dimension $ n \geq 5 $), see [5,7].

\bigskip

Our main result is:

\bigskip

{\it {\bf Theorem.} For all $ b > 0 $, there is a constant $ c $ which depends only on $ b $, $ c=c(b), $ such that,

$$ \sup_{{\mathbb S}_4} u + \inf_{{\mathbb S}_4} u \geq c,  $$

for all solutions $ u $ of $ (E) $.}

\bigskip

\section{Proof of the Theorem.} 

\bigskip

We multiply the equation by $ u $ and integrate:

\be \int_{{\mathbb S}_4} (\Delta u)^2 + 2 \int_{{\mathbb S}_4} |\nabla u|^2 + 6 \int_{{\mathbb S}_4} u = \int_{{\mathbb S}_4} u Ve^{4u} \leq 6 |{\mathbb S}_4| \sup_{{\mathbb S}_4} u. \ee

Let's consider $ G_1 $ and $ G_2 $ the following two Green functions:

$$ \Delta G_1+ 2 G_1= \delta, \,\,\, \Delta G_2 = \delta - \dfrac{1}{|{\mathbb S}_4|} . $$

$ G_1 $ and $ G_2 $ are symmetric.

We refer to [1,6] for the construction of the Green function of 2nd order operators.

We set:

$$ G = G_1 * G_2= \int_{{\mathbb S}_4} G_1(x,z)G_2(z,y) dV_g(z). $$

It is easy to see that:

$$ \Delta^2 G + 2 \Delta G = \delta -\dfrac{1}{|{\mathbb S}_4|}. $$

It is clear that, $ G_1\geq 0 $. We choose $ G_2 $ such that:

$$ G_2 \geq 0,\,\,\,   \int_{{\mathbb S}_4} G_2 \equiv cte>0 $$

We write:

\be \inf_{{\mathbb S}_4} u = u(x) = \dfrac{1}{|{\mathbb S}_4|} \int_{{\mathbb S}_4} u + \int_{{\mathbb S}_4} ( V e^{4u} G - 6 G ) \geq \dfrac{1}{|{\mathbb S}_4|} \int_{{\mathbb S}_4} u -C. \ee

We use $ (1) $ and $ (2) $ to obtain:

\be \sup_{{\mathbb S}_4} u + \inf_{{\mathbb S}_4} u \geq \dfrac{2}{|{\mathbb S}_4|} \int_{{\mathbb S}_4} u + \dfrac{2}{ 6 |{\mathbb S}_4|} \int_{{\mathbb S}_4} |\nabla u|^2 +
\dfrac{1}{6 |{\mathbb S}_4|} \int_{{\mathbb S}_4} (\Delta u)^2 - C. \ee

We use the biharmonic version of the Moser-Trudinger inequality:

\be \log \left [\dfrac{1}{|{\mathbb S}_4|} \int_{{\mathbb S}_4} e^{4u} \right ] \leq \dfrac{1}{ 3 |{\mathbb S}_4|} \int_{{\mathbb S}_4} (\Delta u)^2 + \dfrac{2}{ 3 |{\mathbb S}_4|} \int_{{\mathbb S}_4} |\nabla u|^2 + \dfrac{4}{|{\mathbb S}_4|} \int_{{\mathbb S}_4} u.\ee

We combine $ (3) $ and $ (4) $ to obtain:

$$ \sup_{{\mathbb S}_4} u + \inf_{{\mathbb S}_4} u \geq \dfrac{1}{2} \log \left [ \dfrac{1}{|{\mathbb S}_4|} \int_{{\mathbb S}_4} e^{4u} \right ] - C, $$ 

We integrate the eqaution $ (E) $, we obtain:

$$ \dfrac{1}{|{\mathbb S}_4|} \int_{{\mathbb S}_4} V e^{4u}= 6, $$

thus,

$$ \dfrac{1}{|{\mathbb S}_4|} \int_{{\mathbb S}_4} e^{4u}\geq \dfrac{6}{b}. $$

hence,

\be \sup_{{\mathbb S}_4} u + \inf_{{\mathbb S}_4} u \geq \dfrac{1}{2} \log \left ( \dfrac{6}{b} \right ) - C = \dfrac{\log 6}{2}-\dfrac{\log b}{2}-C. \ee

Remark that:

$$ C \equiv 6\int_{{\mathbb S}_4} G(x,y) dV_g(y) =6 \int_{{\mathbb S}_4} \int_{{\mathbb S}_4} G_1(x,z)G_2(z,y)dV_g(z)dV_g(y)=3 \int_{{\mathbb S}_4} G_2(x,y)dV_g(y). $$

because, $ \int_{{\mathbb S}_4} G_1(x,y) \equiv \dfrac{1}{2}. $

\end{document}